\documentclass[a4paper]{amsart}
\usepackage{bm,color,IEEEtrantools,amsmath,amsthm,amssymb}
\usepackage{setspace}
\usepackage[style=authoryear,natbib=true,maxnames=99,hyperref=true,firstinits=true]{biblatex}
\definecolor{darkblue}{rgb}{0,0,.5} 
\usepackage[colorlinks=true,linkcolor=darkblue,pagecolor=darkblue,citecolor=darkblue,filecolor=darkblue,urlcolor=darkblue,menucolor=darkblue]{hyperref}
\bibliography{/home/jochen/git/Dissertation/library}
        \AtEveryBibitem{%
          \clearfield{day}%
          \clearfield{month}%
          \clearfield{endday}%
          \clearfield{endmonth}%
          \clearfield{doi}
          \clearfield{issn}
        }
\newtheorem{theorem}{Theorem}[section]

\theoremstyle{definition}

\newtheorem{example}[theorem]{Example}

\theoremstyle{remark}

\newtheorem{proposition}[theorem]{Proposition}

\subjclass[2010]{Primary 42A82, 42A16}
\keywords{Gegenbauer expansions, Schoenberg coefficients, isotropic positive definite functions on spheres}
\title{From Fourier to Gegenbauer: Relating Schoenberg coefficients in Gegenbauer expansions on spheres}

\author{Jochen Fiedler\\
    Institute of Applied Mathematics\\
    University of Heidelberg}
\thanks{The author's research was supported by the German Research Foundation (DFG) within the programme ``Spatio/Temporal Probabilistic Graphical Models and Applications in Image Analysis'', grant GRK 1653.}
\date{\today}

\begin{document}

\doublespacing
\maketitle
\newcommand{\hsp}{\hspace{3mm}}
\newcommand{\R}{\mathbb{R}}
\newcommand{\real}{\mathbb{R}}
\newcommand{\mS}{\mathbb{S}}
\newcommand{\N}{\mathbb{N}}
\newcommand{\E}{\mathbb{E}}
\newcommand{\C}{\mathbb{C}}
\newcommand{\Z}{\mathbb{Z}}
\newcommand{\x}{\mathbf{x}}
\newcommand{\y}{\mathbf{y}}
\newcommand{\X}{\mathbf{X}}
\newcommand{\Y}{\mathbf{Y}}
\newcommand{\bZ}{\mathbf{Z}}
\newcommand{\Lc}{\mathcal{L}}
\newcommand{\my}{\mathbf{\mu}}
\newcommand{\sX}{\mathscr{X}}
\newcommand{\sY}{\mathscr{Y}}
\newcommand{\sZ}{\mathscr{Z}}
\newcommand{\Sig}{\mathbf{\Sigma}}
\newcommand{\cov}{\text{{\rm Cov}}}
\newcommand{\var}{\text{Var}}
\newcommand{\corr}{\text{Corr}}
\newcommand{\lo}{L^1(\R^n)}
\newcommand{\lt}{L^2(\R^n)}
\newcommand{\id}{\mathbbmss{1}}
\newcommand{\dit}{\, \text{{\rm d}}}
\newcommand{\eg}{\mbox{e.\,g.}\xspace}
\newcommand{\se}{\text{se}}

\begin{abstract}
It is well-known that every continuous function $\xi:[0,\pi]\rightarrow \real$ admits a series expansion in terms of Gegenbauer polynomials $C_n^{(d-1)/2}$ with coefficients $b_{n,d}$, which are the so-called {\em $d$-dimensional Schoenberg coefficients}. Schoenberg coefficients play an important role in the theory of isotropic positive definite functions on $\mS^d$, since positive definiteness can be characterized by the nonnegativity of the $d$-dimensional Schoenberg coefficients.

In this article we present relations between Schoenberg coefficients of different dimensions. Specifically, we show that the even- resp. odd-dimensional Schoenberg coefficients can be expressed as linear combinations of $1$- resp. $2$-dimensional Schoenberg coefficients, and we give closed form expressions for the coefficients involved in these expansions.
\end{abstract}

\section{Introduction}
Every continuous real-valued function on the interval $[-1,1]$ can be expressed as an infinite series consisting on Gegenbauer polynomials $C_n^\lambda,\;\lambda>-1,n\geq 0$, see for example \citet{Szegö1959}, chapter $3$. Hence, any continuous function $\xi:[0,\pi]\rightarrow \real$ satisfying $\xi(0)=1$ admits for every integer $d\geq 1$ the following {\em $d$-Gegenbauer expansion}
\begin{align}
  \label{eq:gegenbauer}
\xi(\theta)=\sum_{n=0}^\infty b_{n,d}\frac{C_n^{(d-1)/2}(\cos(\theta))}{C_n^{(d-1)/2}(1)},\qquad \theta\in[0,\pi],
\end{align}
where $b_{n,d}$ are called the {\em $d$-dimensional Schoenberg coefficients of $\xi$}. By $\Xi_d$ we denote the class of all those functions $\xi$ for which the $d$-dimensional Schoenberg coefficients are absolutely summable, i.\,e. $\sum_{n=0}^\infty |b_{n,d}|<\infty$.

The motivation of studying $d$-Gegenbauer expansions and their $d$-dimensional Schoenberg coefficients comes from the theory of isotropic positive definite functions on spheres, as we will explain in the following.

For an integer $d\geq 1$ we denote the unit sphere in Euclidean space $\R^{d+1}$ equipped with the Euclidean norm by $\mS^{d}:=\{x\in\R^{d+1}:\lVert x\rVert =1\}$. Consider a kernel $h~:~\mS^d\times \mS^d\rightarrow \real$, which is said to be {\em isotropic} if there exists a function $\xi~:~[0,\pi]\rightarrow \R$ such that
\[
h(x,y)=\xi(\theta(x,y)),\qquad x,y\in\mS^{d},
\]
where $\theta(x,y)=\arccos(\langle x,y\rangle)$ denotes the great circle distance between $x$ and $y$ and $\langle.,.\rangle$ the standard scalar product in $\R^{d+1}$.
The kernel $h:\mS^{d}\times \mS^{d}\rightarrow \R$ is {\em positive definite} if
\begin{align}
  \label{eq:pd}
  \sum_{i=1}^n\sum_{j=1}^na_ia_jh(x_i,x_j)\geq 0,
\end{align}
for all integers $n\geq 1$ and for every choice of constants $a_1,\dots,a_n\in\R$ and every choice of pairwise distinct points $x_1,\dots,x_n\in\mS^{d}$. If the inequality in (\ref{eq:pd}) is strict we call the function $h$ {\em strictly positive definite}. 

We refer to $\Psi_d\; (\Psi_d^+), d=1,2,\dots,$ as the class of continuous functions $\psi:[0,\pi]\rightarrow \R$ with $\psi(0)=1$ for which the associated isotropic kernel $h(x,y)=\psi(\theta(x,y))$ is positive definite (strictly positive definite).

Isotropic positive definite functions on spheres have attracted interest in several areas. They occur as correlation functions for stationary and isotropic random fields on the sphere (\cite{Jones1963}) and, hence, have been studied in spatial statistics (\cite{Banerjee2005,Huang2011} or \cite{Hansen2011}).
Furthermore, they are used  as radial basis functions for interpolating scattered data on spherical domains, see for example \citet{Xu1992, Fasshauer1998} or \citet{Cavoretto2010}. 
Recently, \citet{Gneiting2013, Gneiting2013a} has reviewed conditions for functions to belong to $\Psi_d$ or $\Psi_d^+$ and used them to study parametric families of isotropic and stationary correlation functions on spheres. In his work, he also states several problems for future research, one of which has been solved in \citet{Ziegel2013}, and the solution to another is given here.

Members of $\Psi_d$ and $\Psi^+_d$ are characterized by their $d$-Gegenbauer expansion, see \citet{Schoenberg1942} and \citet{Chen2003}. In particular, the class $\Psi_d,\; d\geq 1,$ consists of functions of form (\ref{eq:gegenbauer}) with $b_{n,d}\geq 0$ and $\sum_{n=0}^\infty b_{n,d}=1$ (this implies $\Psi_d\subset \Xi_d$). For $d\geq 2$, the class $\Psi_d^+$ consists of those functions in $\Psi_d$ for which $b_{n,d}>0$ for infinitely many even and infinitely many odd integers $n$. Consequently, it is possible to study properties of the members of $\Psi_d$ or $\Psi_d^+$ via the coefficients $b_{n,d}$ of the $d$-Gegenbauer expansion (\ref{eq:gegenbauer}). For example, in the cases $d=1$ (\cite{Lorentz1948}) and $d=2$ (\cite{Lang2013}) it was shown that H\"older continuity and differentiability of a function in $\Psi_d$ is connected to the decay rate of $b_{n,d}$.

Since
\[
\Psi_1\supset \Psi_2\supset \Psi_3\supset \dots,
\]
every function in $\Psi_d,\; d\geq 1,$ allows a $1$-Gegenbauer expansion in terms of $C_n^0(\cos\theta)=\cos(n\theta),\;n=0,1,2,\dots$ and coefficients $b_{n,1}$, which is a Fourier cosine expansion. Similarly, since $C^{1/2}_n=P_n$ is a Legendre polynomial, every function $\psi\in\Psi_d$ allows an expansion in terms of Legendre polynomials if $d\geq 2$. Hence, it is interesting to ask how to express higher dimensional Schoenberg coefficients in terms of Fourier or Legendre coefficients.

In general, the connections between Schoenberg coefficients of different dimensions can be helpful to decide, whether or not a function $\psi$ belongs to $\Psi_d$ or $\Psi_d^+$ for a certain $d$.

The problem of expressing even and odd dimensional Schoenberg coefficients in terms of Fourier and Legendre coefficients can be answered using the following recursive identities, stated as Corollary 3 in \citet{Gneiting2013}.\footnote{Note that its proof does not require the Schoenberg coefficients belong to the $d$-Gegenbauer expansion of a positive definite kernel and, hence, it holds also for the Schoenberg coefficients corresponding to $d$-Gegenbauer expansions of members of $\Xi_d$.}
This result provides a connection between $d$-dimensional Schoenberg coefficients and lower dimensional ones. In particular, for all integers $n\geq 1$ it is true that
\begin{align}\label{eq:recursive}
     b_{0,3}=b_{0,1}-\frac12b_{2,1}\qquad\text{and}\qquad   b_{n,3}=\frac12(n+1)(b_{n,1}-b_{n+2,1}).
  \end{align}
Furthermore, if $d\geq 2$, then for all integers $n\geq 0$
\begin{align}
  \label{eq:recursive2}
  b_{n,d+2}=\frac{(n+d-1)(n+d)}{d(2n+d-1)}b_{n,d}-\frac{(n+1)(n+2)}{d(2n+d+3)}b_{n+2,d}
\end{align}
These recursive relationships show that it is possible to express $b_{n,2k+1},k\geq 1,$ as a linear combination of Fourier coefficients $b_{n,1},b_{n+2,1},\dots,b_{n+2k,1}$. Similarly, we can express $b_{n,2k+2},k\geq1,$ as a linear combination of Legendre coefficients $b_{n,2},b_{n+2,2},\dots,b_{n+2k,2}$.

The aim of this work is to provide closed form expressions of the coefficients appearing in these linear combinations, which was stated as Problem 1 in \citet{Gneiting2013a}.

\section{Main results}
In this section we give explicit expressions for Schoenberg coefficients in terms of Fourier cosine and Legendre coefficients. The proofs are provided in Sections \ref{sect:proof1} and \ref{sect:proof2}, respectively.

\begin{theorem}
  \label{thm:main}
  For integers $k\geq 1$ and $n\geq 0$ the Schoenberg coefficient $b_{n,2k+1}$ of the $(2k-1)$-Gegenbauer expansion of a function $\xi\in\Xi_d$ can be expressed in terms of its Fourier cosine coefficients $b_{n,1},b_{n+2,1},\dots,b_{n+2k,1}$, in that
  \[
  b_{n,2k+1}=\sum_{i=0}^ka_i(n,k)b_{n+2i,1},
  \]
  where the $a_i(n,k)$ are given by 
  \begin{align}
    \label{eq:main}
    a_i(n,k)=\frac{(-1)^i}{2^k}{k\choose i}\frac{(n+k)(n+2i)}{(2k-1)!!}\frac{(n+1)_{(2k-1)}}{(n+i)_{(k+1)}},
    \end{align}
    for $(i,n)\neq (0,0)$, whereas $a_0(0,k)=1$ if $i=n=0$. Here $(2k-1)!!=\prod_{i=1}^k(2i-1)$ denotes the {\em double factorial} and $(x)_{(m)}=x(x+1)\cdots(x+m-1)$ the {\em Pochhammer symbol}.
\end{theorem}

\begin{example}
  Consider $k=4$. For $n> 0$ we get 
\begin{align*}
  a_0(n,4)&=\kappa(n+4)(n+5)(n+6)(n+7)\\
  a_1(n,4)&=-4\kappa(n+2)(n+4)(n+6)(n+7)\\
  a_2(n,4)&=6\kappa(n+1)(n+4)^2(n+7)\\
  a_3(n,4)&=-4\kappa(n+1)(n+2)(n+4)(n+6),\\
  a_4(n,4)&=\kappa(n+1)(n+2)(n+3)(n+4),
\end{align*}
where $\kappa=\frac1{1680}$. 
\end{example}

One sees that $a_0(n,4)$ and $a_4(n,4)$ can be expressed in even simpler forms. In general, for $i=0$ and $i=k$ equation (\ref{eq:main}) reduces to
   \begin{align*}
     a_0(n,k)&=\frac{1}{2^k(2k-1)!!}(n+k)_{(k)}\\
     \text{and}\hspace{4cm}&\\
    a_k(n,k)&=\left(-\frac12\right)^k\frac1{(2k-1)!!}(n+1)_{(k)},
  \end{align*}
  respectively.
  
  It is interesting to note that the value of $\sum_{i=0}^k a_i(n,k)$ is either $0$ or $\frac12$.
  \begin{proposition}
    \label{rem:sum1}
    For all integers $k\geq1$ it is true that
    \begin{align*}
    \sum_{i=0}^k a_i(n,k)=\begin{cases}0, & n>0,\\ \frac12, & n=0.\end{cases}
    \end{align*}
    The proof can be found in Section \ref{sect:proof1}.
  \end{proposition}
  
Now let us turn to the analogous problem of finding an expression for $b_{n,2k+2}, k\geq 1,$ in terms of the Legendre coefficients $b_{n,2},\dots,b_{n+2k,2}$.  
\begin{theorem}
  \label{thm:main2}
  For all integers $k\geq 1$ and $n\geq 0$ it is true that
  \begin{align}
    b_{n,2k+2}=\sum_{i=0}^ku_i(n,k)b_{n+2i,2},
  \end{align}
  where $u_i(n,k)$ are given by
  \begin{align}
    \label{eq:main2}
  u_i(n,k)=(-1)^i\frac{(2k-1)!!}{2^k}{k\choose i}{2k+n\choose n}\frac1{(n+i+1/2)_{(k-i)}(n+k+3/2)_{(i)}}.
\end{align} 
\end{theorem}
The proof is provided in Section \ref{sect:proof2}.

\subsection{Applications}
The following example shows that our results can be used to decide whether a function $\psi\in\Psi_1$  is a member of $\Psi_\infty^+$. 

\begin{example}
  Let
  \[
  b_{n,1}=\frac3{\pi^2n^2}, \qquad n\geq 1,
  \]
  and $b_{0,1}=\frac12$. Evidently, the corresponding function $\psi$ is in $\Psi_1$. 
  A symbolical calculation with Mathematica yields for $n,k\geq 1$
  \begin{align*}
    b_{n,2k+1}=\sum_{i=0}^ka_i(n,k)b_{n+2i,1}=\frac{3k(n+k)B(n/2,k)^2}{2n\pi^2(n+2k)^2B(n,2k)},
  \end{align*}
  where $B(x,y)$ denotes the Beta function, and for $n=0,k\geq 1$ it yields
  \begin{align*}
    b_{0,2k+1}=\frac{2(k+1)\pi^2-3k\; _4F_3(1,1,1,1-k;\,2,2,2+k;\,1)}{4 (1 + k) \pi^2},
  \end{align*}
  where $_4F_3$ denotes a generalized hypergeometric function, see \citet{Slater1966}.
  We see that $b_{n,2k+1}>0$ for all $n,k\geq 1$. Now by the definition of the generalized hypergeometric function it is
  \begin{align*}
    _4F_3(1,1,1,1-k;\,2,2,2+k;\,1)&=\sum_{i=0}^\infty\frac{(1)_i(1)_i(1)_i(1-k)_i}{(2)_i(2)_i(2+k)_i}\frac1{i!}\\
    &=\sum_{i=0}^{k-1}\frac1{(1+i)^2}\frac{(1-k)_i}{(2+k)_i},
  \end{align*}
  because $\frac{(1-k)_i}{(2+k)_i}=0$ for $i\geq k$. Furthermore, we have $\frac{(1-k)_i}{(2+k)_i}\leq 1$ for all integers $i,k\geq 0$, and this gives us
  \begin{align*}
    \sum_{i=0}^{k-1}\frac1{(1+i)^2}\frac{(1-k)_i}{(2+k)_i}\leq \sum_{i=0}^{k-1}\frac1{(1+i)^2}\leq \sum_{i=0}^{\infty}\frac1{(1+i)^2}=\frac{\pi^2}{6}.
  \end{align*}
  Hence, we see that
  \begin{align*}
    b_{0,2k+1}&=\frac{2(k+1)\pi^2-3k\; _4F_3(1,1,1,1-k;\,2,2,2+k;\,1)}{4 (1 + k) \pi^2}\\
    &\geq \frac{2(k+1)\pi^2-k\pi^2/2}{4 (1 + k) \pi^2}=\frac{3k+4}{8 (1 + k)}>0,
  \end{align*}
  for all $k\geq 0$. Consequently, $\psi\in\Psi_\infty^+$.
  
  It is interesting to note that $b_{n,2k+1}=\mathcal O(n^{-2})$ for all $k\geq 1$ and, hence, the Schoenberg coefficients show the same asymptotic behaviour in every odd dimension, which can be seen as follows. Stirling's formula for Gamma functions (see 6.1.37 in \cite{Abramowitz1972}) yields for fixed $y$
  \begin{align*}
    B(x,y)&=\Gamma(y)e^y\left(1+\frac{y}x\right)^{1/2-y-x}x^{-y}\frac{1+\mathcal O(x^{-1})}{1+\mathcal O((x+y)^{-1})},
  \end{align*}
  where we used the well-known {\em big $\mathcal O$ notation}.
  Because $\left(1+\frac{y}x\right)^{1/2-y-x}\rightarrow e^{-y}>0$ if $x\rightarrow \infty$, it follows that
  \begin{align*}
    B(x,y)=\Gamma(y)x^{-y}\mathcal O(1)\frac{1+\mathcal O(x^{-1})}{1+\mathcal O((x+y)^{-1})},
  \end{align*}
  yielding immediately
  \[
  \frac{B(n/2,k)^2}{B(n,2k)}=\mathcal O(1)
  \]
  and consequently $b_{n,2k+1}=\mathcal O(n^{-2})$.
\end{example}

\section{Proof of Theorem \ref{thm:main} and Proposition \ref{rem:sum1}}
\label{sect:proof1}
\begin{proof}[{\bf Proof of Theorem \ref{thm:main}:}]
  First we consider the case $n>0$. We proceed by induction over $k\geq 1$.

  Let $k=1$. For all $n\geq 1$ we have
  \[
  b_{n,3}=\frac12(n+1)(b_{n,1}-b_{n+2,1}),
  \]
  yielding $a_0(n,1)=\frac12(n+1)$ and $a_1(n,1)=-\frac12(n+1)$. Inserting $i=0$ and $k=1$ into formula (\ref{eq:main}) immediately yields the same results, proving the claim for $k=1$.

Suppose we have proven (\ref{eq:main}) for an arbitrary $k\geq 1$. From this we want to deduce (\ref{eq:main}) for $k+1$. With (\ref{eq:recursive2}) we see, by comparing coefficients, that 
\begin{IEEEeqnarray*}{rCl}
  \lefteqn{b_{n,2(k+1)+1}=\frac{(n+2k)_{(2)}}{2(2k+1)(n+k)}b_{n,2k+1}-\frac{(n+1)_{(2)}}{2(2k+1)(n+k+2)}b_{n+2,2k+1}}\\
  &=&\frac{(n+2k)_{(2)}}{2(2k+1)(n+k)}\sum_{i=0}^ka_i(n,k)b_{n+2i,1}-\frac{(n+1)_{(2)}}{2(2k+1)(n+k+2)}\sum_{i=0}^ka_i(n+2,k)b_{n+2+2i,1}\\
  &=&\frac{(n+2k)_{(2)}}{2(2k+1)(n+k)}a_0(n,k)b_{n,1}-\frac{(n+1)_{(2)}}{2(2k+1)(n+k+2)}a_k(n+2,k)b_{n+2(k+1),1}\\
&&+\frac1{2(2k+1)}\sum_{i=1}^{k}\left[\frac{(n+2k)_{(2)}}{n+k}a_i(n,k)-\frac{(n+1)_{(2)}}{n+k+2}a_{i-1}(n+2,k)\right]b_{n+2i,1}.
\end{IEEEeqnarray*}
Using the induction hypothesis and the trivial identity 
\begin{align}\label{eq:pochhammer}
  (x)_{(k)}(x+k)_{(l)}=(x)_{(k+l)}\qquad\text{for integers }k,l\geq 0,
\end{align}
we see that
\begin{align*}
  \frac{(n+2k)_{(2)}}{2(2k+1)(n+k)}a_0(n,k)&=\frac{(n+2k)_{(2)}}{2(2k+1)(n+k)}\frac1{2^k}\frac{(n+k)n(n+1)_{(2k-1)}}{(2k-1)!!(n)_{(k+1)}}\\
  &=\frac1{2^{k+1}}\frac{n(n+k+1)(n+1)_{(2k+1)}}{(2k+1)!!(n)_{(k+2)}},
\end{align*}
proving the validity of (\ref{eq:main}) for $i=0$, and
\begin{IEEEeqnarray*}{rCl}
  \lefteqn{-\frac{(n+1)_{(2)}}{2(2k+1)(n+k+2)}a_k(n+2,k)}\\
  &=&-\frac{(n+1)_{(2)}}{2(2k+1)(n+k+2)}\frac{(-1)^k}{2^k}\frac{(n+2+k)(n+2+2k)(n+2+1)_{(2k-1)}}{(2k-1)!!(n+2+k)_{(k+1)}}\\
  &=&\frac{(-1)^{k+1}}{2^{k+1}}\frac{(n+2+2k)(n+1)_{(2k+1)}}{(2k+1)!!(n+2+k)_{(k+1)}}\\
  &=&\frac{(-1)^{k+1}}{2^{k+1}}\frac{(n+k+1)(n+2k+2)(n+1)_{(2k+1)}}{(2k+1)!!(n+1+k)_{(k+2)}},
\end{IEEEeqnarray*}
which shows the validity of (\ref{eq:main}) for $i=k+1$.

It remains to show for $1\leq i\leq k$ that 
\begin{align}
  \label{eq:2}
  \begin{aligned}
    &\frac{(-1)^i}{2^{k+1}}{k+1\choose i}\frac{(n+k+1)(n+2i)(n+1)_{(2k+1)}}{(2k+1)!!(n+i)_{(k+2)}}\\
    &=\frac1{2(2k+1)}\left[\frac{(n+2k)_{(2)}}{n+k}a_i(n,k)-\frac{(n+1)_{(2)}}{n+k+2}a_{i-1}(n+2,k)\right],
  \end{aligned}
\end{align}
where $a_i(n,k)$ and $a_{i-1}(n+2,k)$ can be expressed as in (\ref{eq:main}).
Plugging the induction hypothesis into (\ref{eq:2}) we can reformulate this as
\begin{IEEEeqnarray*}{rCl}
\lefteqn{\frac{(-1)^i}{2^{k+1}}{k+1\choose i}\frac{(n+k+1)(n+2i)(n+1)_{(2k+1)}}{(2k+1)!!(n+i)_{(k+2)}}}\\
&=&\frac{(-1)^i}{2^{k+1}}\frac{(n+2k)_{(2)}}{n+k}{k\choose i}\frac{(n+k)(n+2i)(n+1)_{(2k-1)}}{(2k+1)!!(n+i)_{(k+1)}}\\
&&-\frac{(-1)^{i-1}}{2^{k+1}}\frac{(n+1)_{(2)}}{n+k+2}{k\choose i-1}\frac{(n+2+k)(n+2i)(n+3)_{(2k-1)}}{(2k+1)!!(n+1+i)_{(k+1)}}.
\end{IEEEeqnarray*}
By using (\ref{eq:pochhammer}) and canceling factors we see that this is equivalent to
\begin{multline*}
{k+1\choose i}(n+k+1)\frac1{(n+i)_{(k+2)}}={k\choose i}\frac1{(n+i)_{(k+1)}}+{k\choose i-1}\frac1{(n+i+1)_{(k+1)}}.
\end{multline*}
After multiplying  with $(n+i)_{(k+2)}$ it remains to show that
\begin{equation}\label{eq:5}
{k+1\choose i}(n+k+1)={k\choose i}(n+i+k+1)+{k\choose i-1}(n+i).
\end{equation}
The right hand side of (\ref{eq:5}) equals
\begin{align*}
  &\frac{k!}{(k-i)!(i-1)!}\left(\frac{n+i+k+1}i+\frac{n+i}{k-i+1}\right)\\
   &=\frac{k!}{(k+1-i)!i!}\big[(k-i+1)(n+i+k+1)+i(n+i)\big]\\
  &=\frac{k!}{(k+1-i)!i!}(k+1) (n+k+1)={k+1\choose i}(n+k+1),
\end{align*}
showing the validity of (\ref{eq:5}) and we are done.

Now we turn to the case $n=0$. Equations (\ref{eq:recursive}) and (\ref{eq:recursive2}) show that $a_i(0,k)$ is given by (\ref{eq:main}) for all $i>0$. To find $a_0(0,k)$, note that for $n=0$ equations (\ref{eq:recursive2}) together with (\ref{eq:recursive}) yield 
 \begin{align*}
   b_{0,d+2}=b_{0,d}-\frac2{d(d+3)}b_{2,d},
 \end{align*}
  for all $d\geq 1$.
  Using (\ref{eq:main}) for $b_{2,d}$, we see that for all $d\geq 1$ it holds that
 \[
 b_{0,d+2}=b_{0,1}-R,
 \]
 where the remainder term $R$ does not depend on $b_{0,1}$. This shows that $a_0(0,k)=1$ for all $k\geq 1$. 
\end{proof}

\begin{proof}[{\bf Proof of Proposition \ref{rem:sum1}:}]
 Let $n>0$. Using $(n+i)_{(k+1)}={n+i+k\choose k+1}(k+1)!$ and Theorem \ref{thm:main} we see that
  \begin{align*}
    \sum_{i=0}^ka_i(n,k)&=\frac{(n+k)(n+1)_{(2k-1)}}{2^k(2k-1)!!}\sum_{i=0}^k(-1)^i{k\choose i}\frac{n+2i}{(n+i)_{(k+1)}}\\
    &=\frac{(n+k)(n+1)_{(2k-1)}}{2^k(k+1)!(2k-1)!!}\sum_{i=0}^k(-1)^i{k\choose i}\frac{n+2i}{{n+i+k\choose k+1}}.
  \end{align*}
  Hence, it suffices to prove that
  \begin{align*}
    \sum_{i=0}^k(-1)^i{k\choose i}\frac{n+2i}{{n+i+k\choose k+1}}=0,
  \end{align*}
  which is equivalent to
  \begin{align}\label{eq:prrem2}
    \sum_{i=0}^k(-1)^i{k\choose i}\frac{n}{{n+i+k\choose k+1}}=-2\sum_{i=0}^k(-1)^i{k\choose i}\frac{i}{{n+i+k\choose k+1}}.
  \end{align}
  Now the left-hand side of (\ref{eq:prrem2}) equals
  \begin{align*}
    \sum_{i=0}^k(-1)^i{k\choose i}\frac{n}{{n+i+k\choose k+1}}&=n\frac{k+1}{2k+1}\frac1{{2k+n\choose n-1}}=\frac{k+1}{{2k+n\choose n}},
  \end{align*}
  where the first equality is due to the following result of R. Frisch which can be found, for example, as Note 21 in \citet{Netto1927}, in that
  \begin{align}\label{eq:frisch}
    \sum_{i=0}^k(-1)^i{k\choose i}\frac1{{b+i\choose c}}=\frac{c}{k+c}\frac1{{k+b\choose b-c}},
  \end{align}
  where $b\geq c$ are positive integers.
  
  For the right-hand side of (\ref{eq:prrem2}) we get in a very similar way
  \begin{align*}
    -2\sum_{i=0}^k(-1)^i{k\choose i}\frac{i}{{n+i+k\choose k+1}}&=-2\sum_{i=0}^k(-1)^i\frac{k!}{(k-i)!(i-1)!}\frac1{{n+i+k\choose k+1}}\\
    &=-2k\sum_{i=0}^k(-1)^i{k-1\choose i-1}\frac1{{n+i+k\choose k+1}}\\
    &=-2k\sum_{i=1}^k(-1)^i{k-1\choose i-1}\frac1{{n+i+k\choose k+1}}\\
    &=2k\sum_{i=0}^{k-1}(-1)^i{k-1\choose i}\frac1{{n+i+1+k\choose k+1}}\\
    &=2k\frac{k+1}{2k}\frac1{{2k+n\choose n}}=\frac{k+1}{{2k+n\choose n}},
  \end{align*}
  where we use (\ref{eq:frisch}) for the final equality, thereby showing  (\ref{eq:prrem2}).

  Now consider the case $n=0$. For $i>0$, equation (\ref{eq:main}) simplifies to
  \[
  a_i(0,k)=(-1)^i{k\choose i}\frac1{{k+i\choose k}},
  \]
  which is also valid for $i=0$, since in this case it reduces to $1$. Hence, by using (\ref{eq:frisch}) with $b=c=k$ we get
  \[
  \sum_{i=0}^ka_i(0,k)=\sum_{i=0}^k(-1)^i{k\choose i}\frac1{{k+i\choose k}}=\frac12.
  \]
  \end{proof}

\section{Proof of Theorem \ref{thm:main2}}
\label{sect:proof2}

 We proceed by induction. Let $k=1$. Then
  \[
  b_{n,4}=\frac12(n+1)_{(2)}\left[\frac1{2n+1}b_{n,2}-\frac1{2n+5}b_{n+2,2}\right],
  \]
  implying $u_0(n,1)=\frac1{2(2n+1)}(n+1)_{(2)}$ and $u_1(n,1)=-\frac1{2(2n+5)}(n+1)_{(2)}$. Inserting $k=1,$ and $i=0$ and $i=1$, respectively, in equation (\ref{eq:main2}) proves the claim for $k=1$.

  Suppose we have proven (\ref{eq:main2}) for a $k\geq 1$. We use this to show the validity of (\ref{eq:main2}) for $k+1$. Using (\ref{eq:recursive}) and (\ref{eq:recursive2}) and the induction hypothesis we find that 
\begin{IEEEeqnarray*}{rCl}
 \lefteqn{b_{n,2(k+1)+2} =\frac{(n+2k+1)_{(2)}}{2(k+1)(2n+2k+1)}b_{n,2k+2}-\frac{(n+1)_{(2)}}{2(k+1)(2n+2k+5)}b_{n+2,2k+2}}\\
  &=&\frac{(n+2k+1)_{(2)}}{2(k+1)(2n+2k+1)}\sum_{i=0}^ku_i(n,k)b_{n+2i,2}-\frac{(n+1)_{(2)}}{2(k+1)(2n+2k+5)}\sum_{i=0}^ku_i(n+2,k)b_{n+2+2i,2}\\
  &=&\frac{(n+2k+1)_{(2)}}{2(k+1)(2n+2k+1)}u_0(n,k)b_{n,2}-\frac{(n+1)_{(2)}}{2(k+1)(2n+2k+5)}u_k(n+2,k)b_{n+2(k+1),2}\\
  &&+\sum_{i=1}^k\left[\frac{(n+2k+1)_{(2)}}{2(k+1)(2n+2k+1)}u_i(n,k)-\frac{(n+1)_{(2)}}{2(k+1)(2n+2k+5)}u_{i-1}(n+2,k)\right]b_{n+2i,2}.
\end{IEEEeqnarray*}
Using $(2k-1)!!=\frac{(2k)!}{2^kk!}$ and (\ref{eq:pochhammer}) we see that
\begin{IEEEeqnarray*}{rCl}
  \frac{(n+2k+1)_{(2)}}{2(k+1)(2n+2k+1)}u_0(n,k)&=&\frac{(n+2k+1)_{(2)}(2k-1)!!}{2^{k+1}(k+1)(2n+2k+1)}{2k+n\choose n}\frac1{(n+1/2)_{(k)}}\\
  &=&\frac1{2^{k+1}2^{k+1}}\frac{(2k+2+n)!}{(k+1)!n!}\frac1{(n+1/2)_{(k+1)}}\\
  &=&\frac{(2k+1)!!}{2^{k+1}}{2k+2+n\choose n}\frac1{(n+1/2)_{(k+1)}},
\end{IEEEeqnarray*}
proving the claim for $i=0$. In a very similar way we get
\begin{IEEEeqnarray*}{rCl}
  \lefteqn{-\frac{(n+1)_{(2)}}{2(k+1)(2n+2k+5)}u_k(n+2,k)}\\
  &=&(-1)^{k+1}\frac{(n+1)_{(2)}(2k-1)!!}{2^{k+1}(k+1)(2n+2k+5)}{2k+n+2\choose n+2}\frac1{(n+k+7/2)_{(k)}}\\
   &=&(-1)^{k+1}\frac1{2^{k+1}2^{k+1}}\frac{(2k+2+n)!}{(k+1)!n!}\frac1{(n+k+5/2)_{(k/2)}}\\
  &=&(-1)^{k+1}\frac{(2k+1)!!}{2^{k+1}}{2k+2+n\choose n}\frac1{(n+k+5/2)_{(k/2)}},
\end{IEEEeqnarray*}
confirming the claim for $i=k+1$.

Now let $1\leq i \leq k$. We need to show that
\begin{align*}
  u_i(n,k+1)&=(-1)^i\frac{(2k+1)!!}{2^{k+1}}{k+1\choose i}{2k+2+n\choose n}\frac1{(n+i+1/2)_{(k+1-i)}(n+k+1+3/2)_{(i)}}\\
  &= \frac{(n+2k+1)_{(2)}}{2(k+1)(2n+2k+1)}u_i(n,k)-\frac{(n+1)_{(2)}}{2(k+1)(2n+2k+5)}u_{i-1}(n+2,k),
\end{align*}
which is equivalent to
\begin{align}\label{eq:6}
  \begin{aligned}
   \lefteqn{(-1)^i\frac{(2k+1)!!}{2^{k+1}}{k+1\choose i}{2k+2+n\choose n}\frac1{(n+i+1/2)_{(k+1-i)}(n+k+1+3/2)_{(i)}}}\\
  &=\frac{(-1)^i(n+2k+1)_{(2)}(2k-1)!!}{2^{k+1}(k+1)(2n+2k+1)}{k\choose i}{2k+n\choose n}\frac1{(n+i+1/2)_{(k-i)}(n+k+3/2)_{(i)}}\\
  &-\frac{(-1)^{i-1}(n+1)_{(2)}(2k-1)!!}{2^{k+1}(k+1)(2n+2k+5)}{k\choose i-1}{2k+n+2\choose n+2}\frac1{(n+i+1+1/2)_{(k-i+1)}(n+2+k+3/2)_{(i-1)}}.
\end{aligned}
\end{align}
Using (\ref{eq:pochhammer}) and using similar arguments as in the cases $i=0$ and $i=k+1$ we see that (\ref{eq:6}) is equivalent to
\begin{multline*}
    {k+1\choose i}\frac1{(n+i+1/2)_{(k+1-i)}(n+k+5/2)_{(i)}}\\
  ={k\choose i}\frac1{(n+i+1/2)_{(k+1-i)}(n+k+3/2)_{(i)}}+{k\choose i-1}\frac1{(n+i+3/2)_{(k+1-i)}(n+k+5/2)_{(i)}}.
\end{multline*}
Multiplying with $(n+i+1/2)_{(k+1-i)}(n+k+5/2)_{(i)}$ illustrates that we need to show 
\begin{equation*}
    {k+1\choose i}
  ={k\choose i}\frac{n+k+i+3/2}{n+k+3/2}+{k\choose i-1}\frac{n+i+1/2}{n+k+3/2}.
\end{equation*}
Simplifying the right hand side yields
\begin{multline*}
  {k\choose i}\frac{(n+k+i+3/2)(k-i+1)+i(n+i+1/2)}{(n+k+3/2)(k-i+1)}\\
  ={k\choose i}\frac{(k+1) (n+k+3/2)}{(n+k+3/2)(k-i+1)} ={k+1\choose i},
\end{multline*}
and the proof is complete.
\section*{Acknowledgements}
I would like to thank Tilmann Gneiting for helpful discussions and comments.

\printbibliography

\end{document}